\documentclass[12pt]{amsart}
\usepackage{amssymb,latexsym}

\numberwithin{equation}{section}

\newcommand{\G}{G_{r,p,n}}

\newcommand{\sumlim}{\sum\limits}

\newcommand{\fix}{{\rm fix}}
\newcommand{\csum}{{\rm csum}}
\newcommand{\exc}{{\rm exc}}
\newtheorem{thm}{Theorem}[section]

\newtheorem{lem}[thm]{Lemma}

\newtheorem{cor}[thm]{Corollary}

\newtheorem{obs}[thm]{Observation}

\newtheorem{exa}[thm]{Example}

\theoremstyle{definition}
\newtheorem{defn}[thm]{Definition}

\newtheorem{prop}[thm]{Proposition}
\newtheorem{conj}[thm]{Conjecture}
\newtheorem{clm}[thm]{Claim}
\newcommand{\een}{\end{enumerate}}
\newcommand{\blem}{\begin{lem}}
\newcommand{\elem}{\end{lem}}
\newcommand{\bcl}{\begin{clm}}
\newcommand{\ecl}{\end{clm}}
\newcommand{\bthm}{\begin{thm}}
\newcommand{\ethm}{\end{thm}}
\newcommand{\bpr}{\begin{prop}}
\newcommand{\epr}{\end{prop}}
\newcommand{\bco}{\begin{cor}}
\newcommand{\eco}{\end{cor}}
\newcommand{\bcon}{\begin{conj}}
\newcommand{\econ}{\end{conj}}
\newcommand{\bde}{\begin{defn}}
\newcommand{\ede}{\end{defn}}
\newcommand{\bex}{\begin{exa}}
\newcommand{\eexa}{\end{exa}}
\newcommand{\bobs}{\begin{obs}}
\newcommand{\eobs}{\end{obs}}
\newcommand{\bexe}{\begin{exe}}
\newcommand{\eexe}{\end{exe}}

\newcommand{\grn}{G_{r,n}}

\begin{document}
\pagenumbering{arabic}
\def\qed{\hfill\rule{2mm}{2mm}}

\title[Excedances for involutions in $G_{r,p,n}$]{Excedance number for involutions in complex reflection groups}

\author{Eli Bagno}
\address{Einstein institute of Mathematics, The Hebrew University,
Givat Ram, 91904 Jerusalem, Israel, and The Jerusalem College of
Technology, Jerusalem, Israel}
\email{bagnoe@math.huji.ac.il,bagnoe@jct.ac.il}

\author{David Garber}
\address{School of Sciences, Holon Institute of Technology, PO Box 305,
58102 Holon, Israel} \email{garber@hit.ac.il}

\author{Toufik Mansour}
\address{Department of Mathematics, University of Haifa, Haifa, Israel}
\email{toufik@math.haifa.ac.il}

\date{\today}

\maketitle
\begin{abstract}
We define the excedance number on the complex reflection groups and
compute its multidistribution with the  number of fixed points on
the set of involutions in these groups. We use some recurrence
formulas and generating functions manipulations to obtain our
results.
\end{abstract}

\section{Introduction}

Let $V$ be a complex vector space of dimension $n$.  A {\em
pseudo-reflection} on $V$ is a linear transformation on $V$ of
finite order which fixes a hyperplane in $V$ pointwise. A {\em
complex reflection group} on $V$ is a finite subgroup $W$ of ${\rm
GL}(V)$ generated by pseudo-reflections.

Irreducible finite complex reflection groups have been classified by
Shephard-Todd \cite{ST}. In particular, there is a single infinite
family of groups and exactly 34 other ``exceptional" complex
reflection groups. The infinite family $\G$, where $r,p,n$ are
positive integers with $p | r$, consists of the groups of $n\times
n$ matrices such that:

\begin{enumerate}
\item The entries are either 0 or $r^{\rm th}$ roots of unity;
\item  There is exactly one nonzero entry in each row and each column;
\item  The $(r/p)^{{\rm th}}$ power of the product of the nonzero
entries is 1.
\end{enumerate}

If $p=1$ then we get the colored permutation group:
$G_{r,n}=G_{r,1,n}$. It consists of all permutations of the set
$$\Sigma=\{1,\dots,n,\bar{1},\dots,\bar{n},\dots,1^{[r-1]},\dots,n^{[r-1]}\}$$
satisfying $\pi(\bar{i})=\overline{\pi(i)}$.

The classical Weyl groups appear as special cases: $G_{1,1,n}=S_n$
the symmetric group, $G_{2,1,n}=B_n$ the hyperoctahedral group, and
$G_{2,2,n}=D_n$, the group of even-signed permutations.

In $S_n$ one can define the following well-known parameters: Given
$\sigma \in S_n$, $i \in [n]$ is {\it an excedance of $\sigma$} if
and only if $\sigma(i)>i$. The number of excedances is denoted by
${\rm exc}(\sigma)$. Another natural parameter on $S_n$ is the
number of fixed points, denoted by ${\rm fix}(\sigma)$.

We say that a permutation $\pi \in G_{r,p,n}$ is {\it an involution}
if $\pi^2=1$. Let $I_{r,p,n}$ be the set of involutions in the
complex reflection group $G_{r,p,n}$.

In this paper we are interested in computing the number of
involutions having specific numbers of fixed points and excedances.
We do this by producing recurrence formulas, and computing them
explicitly by the corresponding generating functions.

\medskip
Here are our main results:

\bthm  (See Corollaries \ref{exccolor of grn}, \ref{r odd} and
\ref{r even p devides half r}).

\begin{enumerate}
\item
The number of involutions  $\pi \in G_{r,p,n}$ where $r$ is odd
and $p|r$ with $\exc^{{\rm Clr}}(\pi)=m$ is:

$$\sumlim_{j=\frac{n}{2}}^{n}{(n-j)!\binom{n}{n-j,\ n-j,\ n-2k,\
2k-2n+2j}{(\frac{r}{2})}^{n-j}}$$

\item
The number of involutions  $\pi \in G_{r,p,n}$ where $r$ is even and
$p|\frac{r}{2}$ with $\exc^{{\rm Clr}}(\pi)=m$ is:
$$k!\binom{n}{k,\ k,\ n-2k}{(\frac{r}{2})}^k$$
where $k=\frac{m}{r}$.

\end{enumerate}
\ethm

\bthm (See corollary \ref{p dont divide half r}).

The number of involutions  $\pi \in G_{r,p,n}$ ($r$ is even, $p \not
| \frac{r}{2}$) with $\exc^{{\rm Clr}}(\pi)=m$ is:
\begin{equation*}
\frac{(\frac{m}{r})!}{2^{\frac{m}{r}}}\binom{n}{\frac{m}{r},\
\frac{m}{r},\ n-2\frac{m}{r}}{(r+1)}^{\frac{m}{r}}
\end{equation*}
\ethm

This paper is organized as follows. In Section \ref{pre}, we recall
some properties of $G_{r,p,n}$. In Section \ref{stat} we define some
parameters on $\grn$ and hence also on $G_{r,p,n}$. In Section
\ref{involutions} we classify the involutions of $\grn$ and
$G_{r,p,n}$ and finally in Section \ref{computations} we compute the
corresponding recurrence and explicit formulas.

\section{Preliminaries}\label{pre}

\subsection{Complex reflection groups}\label{grn}

\bde Let $r$ and $n$ be positive integers. {\it The group of colored
permutations of $n$ digits with $r$ colors} is the wreath product
$$\grn=\mathbb{Z}_r \wr S_n=\mathbb{Z}_r^n \rtimes S_n,$$ consisting of all the pairs $(z,\tau)$ where
$z$ is an $n$-tuple of integers between $0$ and $r-1$ and $\tau \in
S_n$. The multiplication is defined by the following rule: For
$z=(z_1,...,z_n)$ and $z'=(z'_1,...,z'_n)$
$$(z,\tau) \cdot (z',\tau')=((z_1+z'_{\tau^{-1}(1)},...,z_n+z'_{\tau^{-1}(n)}),\tau \circ \tau')$$ (here $+$ is
taken modulo $r$). \ede

We use some conventions along this paper. For an element
$\pi=(z,\tau) \in \grn$ with $z=(z_1,...,z_n)$ we write
$z_i(\pi)=z_i$. For $\pi=(z,\tau)$, we denote $|\pi|=(0,\tau), (0
\in \mathbb{Z}_r^n)$. An element $(z,\tau)=((1,0,3,2),(2,1,4,3)) \in
G_{3,4}$ will be written as $(\bar{2} 1 \bar{\bar{\bar{4}}}
\bar{\bar{3}})$.

A much more natural way to present $\grn$ is the following: Consider
the alphabet $\Sigma=\{1,\dots,n,\bar{1},\dots,\bar{n},\dots,
1^{[r-1]},\dots,n^{[r-1]} \}$ as the set $[n]$ colored by the colors
$0,\dots,r-1$. Then, an element of $\grn$ is a {\it colored
permutation}, i.e. a bijection $\pi: \Sigma \rightarrow \Sigma$ such
that $\pi(\bar{i})={\overline{\pi (i)}}$.


For each $p|r$ we define the {\it complex reflection group}:

\begin{equation}\label{def-grpn}
G_{r,p,n}:=\{g \in G_{r,n} \mid \csum(g)\equiv 0  \; {\rm mod} \;
p\}.
\end{equation}
 where $${\rm csum}(\sigma) = \sumlim_{i=1}^n z_i(\sigma).$$

\section{Statistics on $\grn$ and its subgroups}\label{stat}

In this section we define some parameters on $\grn$. $G_{r,p,n}$
inherits all of them. Given any ordered alphabet $\Sigma'$, we
recall the definition of the {\it excedance set} of a permutation
$\pi$ on $\Sigma'$ :
$${\rm Exc}(\pi)=\{i \in \Sigma' \mid \pi(i)>i\}$$ and the {\it excedance
number} is defined to be ${\rm exc}(\pi)=|{ \rm Exc}(\pi)|$.\\

\bde

We define the color order on the set:
$$\Sigma=\{1,\dots,n,\bar{1},\dots,\bar{n},\dots,
1^{[r-1]},\dots,n^{[r-1]} \}$$ by

$$1^{[r-1]} < \cdots < n^{[r-1]} < 1^{[r-2]} < 2^{[r-2]} < \cdots < n^{[r-2]} < \cdots < 1 < \cdots < n.$$
\ede

We note that there are some other possible ways of defining orders
on $\Sigma$, some of them lead to other versions of the excedance
number, see for example \cite{BG}.
\begin{exa}

Given the color order: $$\bar{\bar{1}} < \bar{\bar{2}}
<\bar{\bar{3}} < \bar{1} < \bar{2} <\bar{3} < 1 < 2 < 3,$$ we write
$\sigma=(3\bar{1}\bar{\bar{2}}) \in G_{3,3}$ in an extended form:

$$\begin{pmatrix} \bar{\bar{1}} & \bar{\bar{2}} & \bar{\bar{3}} &
\bar{1} & \bar{2}& \bar{3} & 1 & 2 & 3\\
\bar{\bar{3}} & 1 & \bar{2} & \bar{3} & \bar{\bar{1}}  &  2 & 3 &
\bar{1} & \bar{\bar{2}}
\end{pmatrix}$$
and calculate: ${\rm
Exc}(\sigma)=\{\bar{\bar{1}},\bar{\bar{2}},\bar{\bar{3}},\bar{1},\bar{3},1\}$
and ${\rm exc}(\sigma)=6$.
\end{exa}

Before defining the excedance number, we have to introduce some
notions.

Let $\sigma \in \grn$. We define:
$${\rm csum}(\sigma) = \sumlim_{i=1}^n z_i(\sigma)$$

$${\rm Exc}_A(\sigma) = \{ i \in [n-1] \ | \ \sigma(i) > i \}$$
where the comparison is with respect to the color order.
$${\rm exc}_A(\sigma) = |{\rm Exc}_A(\sigma)|$$

\begin{exa}
Take $\sigma=(\bar{1}\bar{\bar{3}}4\bar{2}) \in G_{3,4}$. Then ${\rm
csum}(\sigma)=4$, ${\rm Exc_A}(\sigma)=\{3\}$ and hence ${\rm
exc}_A(\sigma)=1$.
\end{exa}

Let $\sigma \in \grn$. Recall that for $\sigma=(z,\tau)\in \grn$,
$|\sigma|$ is the permutation of $[n]$ satisfying
$|\sigma|(i)=\tau(i)$. For example, if
$\sigma=(\bar{2}\bar{\bar{3}}1\bar{4})$ then $|\sigma|=(2314)$.

Now we can define the colored excedance number for $\grn$.

\bde Define:

$${\rm exc}^{{\rm Clr}}(\sigma)=r \cdot {\rm exc}_A(\sigma)+{\rm csum}(\sigma)$$
\ede

One can view ${\rm exc}^{{\rm Clr}}(\sigma)$ in a different way (see
\cite{BG}):

\blem Let $\sigma \in \grn$. Consider the set $\Sigma$ ordered by
the color order. Then $${\rm exc}(\sigma)={\rm exc}^{{\rm
Clr}}(\sigma).$$ \elem

We say that $i \in [n]$ is an {\it absolute fixed point} of $\sigma
\in \grn$ if $|\sigma(i)|=i$.

\section{Involutions in $G_{r,p,n}$}\label{involutions}
As was already mentioned, we say that $\sigma$ is an {\it
involution} if $\sigma ^2=1$.

In this section we classify the involutions of $\G$. Note that each
involution of $\G$ can be decomposed into a product of 'atomic'
involutions of two types: absolute fixed points and $2$-cycles.

\medskip

 We start with the absolute fixed points. In the case $p=1$,
i.e. $\mathbb{Z}_r \wr S_n = G_{r,n}$, we split into two subcases
according to the parity of $r$.  In the case of even $r$, an
absolute fixed point can be one of the following two kinds:
$\pi(i)=i$ or $\pi(i)=i^{\left[{r \over 2}\right]}$. If $r$ is odd,
an absolute fixed point can be only of the first kind.

If $p>1$ and $r$ is odd, we have the same absolute fixed points as
in the case $p=1$. On the other hand, if $r$ is even, then we have
to split again into two subcases. If $p|{r \over 2}$, then the
absolute fixed points in $I_{r,p,n}$ are exactly as those of
$I_{r,1,n}$. If $p \hspace{-5pt} \not\hspace{-2pt}| \hspace{2pt} {r
\over 2}$, then an element of with an odd number of fixed points of
the form $\pi(i)=i^{\left[{r \over 2}\right]}$ is not an element of
$I_{r,p,n}$ and thus the only absolute fixed points are of the form
$\pi(i)=i$ or pairs of absolute fixed points of the form:
$\pi(i)=i^{\left[{r \over 2}\right]}; \pi(j)= j^{\left[{r \over
2}\right]}$.

In all cases, the $2$-cycles have the form $\pi(i)= j^{[k]}$;
$\pi(j)= i^{[r-k]}$ where $0 \leq k \leq r-1$.

We conclude this section with an example:

\begin{exa}\label{example of involution}

Let $r=18, p=6, n=7$ and let
$$\pi=\left(\begin{array}{ccccccc} 1 & 2 & 3 & 4 & 5 & 6 & 7 \\
7^{[2]} &3 &2 &4^{[9]} & 5^{[9]}&6 & 1^{[16]}
\end{array}\right) \in I_{18,6,7}.$$
Then $\pi$ can be decomposed into the absolute fixed points: $\left( \begin{array}{c}6 \\
6\end{array}\right)$, the pair of absolute fixed points:
$\left(\begin{array}{cc} 4 & 5 \\4 ^{[9]} & 5^{[9]}
\end{array}\right)$ and the following two 2-cycles: $\left(\begin{array}{cc} 1 & 7 \\
7^{[2]} & 1^{[16]} \end{array}\right)$ and $\left(\begin{array}{cc} 2 & 3 \\
3 & 2 \end{array}\right)$.

\end{exa}

\section{Recurrence and explicit formulas}\label{computations}

In this section, we compute recurrence and explicit formulas for
$$f_{r,p,n}(u,v,w) = \sum_{\pi\in
I_{r,p,n}}u^{\fix(\pi)}v^{\exc(\pi)}w^{\csum(\pi)}$$ for all $r$ and
$p$ where $p|r$.

\subsection{Recurrence formulas for
$G_{r,n}=G_{r,1,n}$}\label{recurrence}
Let $\pi$ be any colored
involution in $I_{r,n}=I_{r,1,n}$. Then we have either
$\pi(n)=n^{[j]}$ or $\pi(n)=k^{[j]}$ with $k<n$.

If $\pi(n)=n^{[j]}$, then we divide into two subcases according to
the parity of $r$, as we have seen in Section \ref{involutions}. If
$r$ is even we have $j=0$ or $j=\frac{r}{2}$. If $r$ is odd then
$j=0$.

For $\pi \in I_{r,1,n}$  such that $\pi(n)=n^{[j]}$, define $\pi'
\in I_{r,1,n-1}$ by ignoring the last digit of $\pi$. For $\pi \in
I_{r,1,n}$ with $\pi(n)=k^{[j]}$ and $\pi(k)=n^{[r-j]}$, define
$\pi'' \in I_{r,1,n-2}$ in the following way: Write $\pi$ in its
complete notation, i.e. as a matrix of two rows, as in Example
\ref{example of involution}. The first row of $\pi''$ is
$(1,2,...,n-2)$ while the second row is obtained from the second row
of $\pi$ by ignoring the digits $n$ and $k$ and the other digits are
placed in an order preserving way with respect to the second row of
$\pi$. Here is an explicit formula for the map $\pi \mapsto \pi''$.

$$\pi'' (i)=\left\{
\begin{array}{ccc}
\pi (i)      & & 1 \leq i < k\ {\quad and \qquad} \pi(i) < k \\
\pi (i)-1    & & 1 \leq i < k\ {\quad and \qquad} \pi(i) > k \\
\pi (i-1)    & &   k \leq i < n \ {\quad and \qquad} \pi(i) < k \\
\pi (i-1)-1  & &   k \leq i < n \ {\quad and \qquad} \pi(i) > k \\
\end{array}\right.$$

Note that the map $\pi \mapsto \pi'$ is a bijection from the set
$$\{\pi \in I_{r,1,n} \mid \pi(n)=n^{[j]}\} \quad (j\ {\rm fixed})$$
to $I_{r,1,n-1}$, while $\pi \mapsto \pi''$ is a bijection from the
set $\{\pi \in I_{r,1,n} \mid \pi(n)=k^{[j]}\} \quad (j\ {\rm
fixed})$ to $I_{r,1,n-2} $.

For any $r$, if $\pi(n)=n^{[j]}$ then:
$$\fix(\pi)=\fix(\pi')+1,$$ $$\exc_A(\pi)=\exc_A(\pi'),$$
$$\csum(\pi)=\csum(\pi')+j. $$

If  $\pi(n)=t^{[j]}$, then the parameters satisfy
$$\fix(\pi)=\fix(\pi''),$$
$$\exc_A(\pi)=\exc_A(\pi'')+\delta_{j,0},$$
$$\csum(\pi)=\csum(\pi'')+r(1-\delta_{j,0}).$$
where $\delta_{i,j}$ is the Kronecker Delta:
$$\delta_{i,j}=
\left\{ \begin{array}{cc} 1 & i=j \\ 0 & i\neq j \end{array}
\right.$$

The above consideration gives the following recurrence formula,
where we define $\mu_{r}=1+w^{r \over 2}$ for even $r$, and
$\mu_{r}=1$ otherwise:

\begin{eqnarray*}
f_{r,1,n}(u,v,w) & = & u\mu_r f_{r,1,n-1}(u,v,w) \\
                    & & +(n-1)(v+(r-1)w^{r})f_{r,1,n-2}(u,v,w),\quad n\geq1
\end{eqnarray*}

\subsection{Explicit formulas for
$G_{r,n}=G_{r,1,n}$}\label{exp_grn}

We turn now to the explicit formula. Define:
$$F_{r,p}(x;u,v,w)=\sumlim_{n \geq 0}{f_{r,p,n}(u,v,w)\frac{x^n}{n!}}=\sum_{n\geq0}\sum_{\pi\in I_{r,p,n}}\left(
u^{\fix(\pi)}v^{\exc_A(\pi)}w^{{\rm
csum}(\pi)}\right)\frac{x^n}{n!}.$$

Rewriting the recurrence formula in terms of generating functions,
we obtain that:

\begin{eqnarray*}
x\frac{\partial}{\partial x}F_{r,1}(x;u,v,w)&=&\sum\limits_{n \geq1}\frac{f_{r,1,n}(u,v,w)}{(n-1)!}x^n=\\
& = & ux \mu_r \sum\limits_{n \geq 1}\frac{x^{n-1}}{(n-1)!}f_{r,1,n-1}(u,v,w)\\
&   & \qquad +x^2(v+(r-1)w^{r})\sum\limits_{n\geq2}\frac{x^{n-2}}{(n-2)!}f_{r,1,n-2}(u,v,w)\\
& = & ux \mu_r F_{r,1}(x;u,v,w)+x^2(v+(r-1)w^{r})F_{r,1}(x;u,v,w)
\end{eqnarray*}

Thus, the generating function $F_{r,1}(x;u,v,w)$ satisfies:
$$\frac{\frac{\partial}{\partial x}F_{r,1}(x;u,v,w)}{F_{r,1}(x;u,v,w)}=u \mu_r +x(v+(r-1)w^{r}).$$

Integrating with respect to $x$ in both sides of the above
differential equation, using the fact that $F_{r,1}(0;u,v,w)=1$, we
obtain the following proposition.

\begin{prop}\label{thmm}
Let $r\geq1$. The generating function $F_{r,1}(x;u,v,w)$ is given by
$$e^{ux \mu_r +\frac{1}{2}x^2(v+(r-1)w^{r})}$$
\end{prop}

We are looking for an explicit expression for the polynomial
$f_{r,1,n}(u,v,w)$. From the definitions we have that
$\frac{f_{r,1,n}(u,v,w)}{n!}$ is the coefficient of $x^n$ in
$F_{r,1}(x;u,v,w)$, namely $[x^n]F_{r,1}(x;u,v,w)$. Computing the
coefficient of $x^n$ in the Maclaurin series of $F_{r,1}(x;u,v,w)$
one gets:
$$\begin{array}{ll}
f_{r,1,n}(u,v,w)
&=n![x^n]\sum_{j\geq0}\left(ux\mu_r+\frac{1}{2}x^2(v+(r-1)w^r)\right)^j/j!\\
&=n![x^n]\sum_{j\geq0}\sum_{i=0}^j\binom{j}{i}\frac{x^{i+j}u^{j-i}(v+(r-1)w^r)^i}{j!2^i}\mu_r^{j-i}\\
&=n!\sumlim_{j=n/2}^n\binom{j}{n-j}\frac{u^{2j-n}(v+(r-1)w^r)^{n-j}}{j!2^{n-j}}\mu_r^{2j-n}
\end{array}$$
Hence, we have the following corollary.
\begin{cor}
The polynomial $f_{r,1,n}(u,v,w)$ is given by

\begin{equation}\label{pol}\sum_{j=n/2}^n (n-j)!\binom{n}{n-j,n-j,2j-n}\frac{u^{2j-n}(v+(r-1)w^r)^{n-j}}{2^{n-j}}\mu_r^{2j-n}.\end{equation}
\end{cor}

If we substitute $w=1$ and compute the coefficient of $u^mv^{\ell}$
in Formula (\ref{pol}), we get the following result.

\begin{cor} Let $r\geq1$.  The number of colored involutions in
$G_{r,n}$ with exactly $m$ absolute fixed points and ${\rm exc}_A
(\pi)=\ell$ is given by
$$(\frac{n-m}{2})!(r-1)^{\frac{n-m}{2}-\ell}\binom{n}{\frac{n-m}{2},m,\frac{n-m}{2}-\ell,\ell}
2^{\frac{m(3-2k)-n}{2}}$$ where $k \in \{0,1\}$ and $k \equiv r
\pmod 2$.
\end{cor}

It is easy to see that if $r=1$, then $2\exc_A(\pi)+\fix(\pi)=n$ for
each involution $\pi$ of $S_n$. From the above corollary we have
then that the number of involutions in $S_n$ with exactly $\ell$
excedances is given by $\frac{l!}{2^l}\binom{n}{l,l,n-2l}$

We turn now to the computation of the number of involutions with a
fixed number of excedances. We do this by substituting $u=1$ and
$v=w^r$ in Formula (\ref{pol}).

\begin{cor}\label{exccolor of grn}

The number of involutions  $\pi \in \grn$ with $\exc^{{\rm
Clr}}(\pi)=m$ is:
\begin{equation*}\left\{\begin{array}{cc}
k!\binom{n}{k,\ k,\ n-2k}{(\frac{r}{2})}^k & \qquad r \equiv 1 \pmod 2 \\
\sumlim_{j=\frac{n}{2}}^{n}{(n-j)!\binom{n}{n-j,\ n-j,\ n-2k,\
2k-2n+2j}{(\frac{r}{2})}^{n-j}} & \qquad r \equiv 0 \pmod 2
\end{array} \right.
\end{equation*}
where $k=\frac{m}{r}$.
\end{cor}

Note that $k$ is an integral number, since $\exc^{\rm Clr}(\pi)$ is
an integral multiplicity of $r$, for  $\pi \in I_{r,n}$.

\subsection{Recurrence and explicit formulas for $G_{r,p,n}$ where $r$ is odd, $p>1$}

As we have seen in Section \ref{involutions}, the involutions in
this case coincide with the involutions of $G_{r,1,n}$ where $r$ is
odd and thus we have:

\begin{cor}\label{r odd}
The recurrence formula for $f_{r,p,n}(u,v,w)$ for odd $r$ is:
\begin{eqnarray*}
f_{r,p,n}(u,v,w) & = & u f_{r,p,n-1}(u,v,w) \\
                    & & +(n-1)(v+(r-1)w^{r})f_{r,p,n-2}(u,v,w),\quad n\geq1
\end{eqnarray*}
and thus its explicit formula is:
$$f_{r,p,n}(u,v,w)=\sum_{j=n/2}^n(n-j)!\binom{n}{n-j,\ n-j, \ 2j-n}\frac{u^{2j-n}(v+(r-1)w^r)^{n-j}}{2^{n-j}}.$$
\end{cor}

\subsection{Recurrence and explicit formulas for $G_{r,p,n}$ where $r$ is even and $p>1, p|\frac{r}{2}$}
Also in this case, we have that the involutions coincide with the
involutions of $G_{r,1,n}$ where $r$ is even, and thus we have:

\begin{cor}\label{r even p devides half r}
The recurrence formula for $f_{r,p,n}(u,v,w)$ for even $r$  and
$p>1, p|\frac{r}{2}$ is:
\begin{eqnarray*}
f_{r,p,n}(u,v,w) & = & u (1+w^{\frac{r}{2}})f_{r,p,n-1}(u,v,w) \\
                    & & +(n-1)(v+(r-1)w^{r})f_{r,p,n-2}(u,v,w),\quad n\geq1
\end{eqnarray*}
and thus its explicit formula is:
\begin{tiny}
$$f_{r,p,n}(u,v,w)=\sum_{j=n/2}^n(n-j)!\binom{n}{n-j,\ n-j,\ 2j-n}\frac{u^{2j-n}(v+(r-1)w^r)^{n-j}}{2^{n-j}}(1+w^{\frac{r}{2}})^{2j-n}.$$
\end{tiny}
\end{cor}

\subsection{Recurrence and explicit formulas for $G_{r,p,n}$ where $r$ is even and $p>1, p\not \hspace{-3pt}|\frac{r}{2}$}
Let $\pi$ be any colored involution in $I_{r,p,n}$. Then, according
to Section \ref{involutions}, we have either $\pi(n)=n^{[j]}$ (where
$j=0$ or $j=\frac{r}{2}$) or $\pi(n)=k^{[j]}$ with $k<n$.

We start with the recurrence formula. Let $\pi$ be any colored
involution in $I_{r,p,n}$. Then we have several cases:
\begin{enumerate}
\item $\pi(n)=n$. In this case define $\pi' \in I_{r,p,n-1}$ by ignoring the last digit of $\pi$.
The map $\pi \mapsto \pi'$ is a bijection from the set $\{\pi \in
I_{r,p,n} \mid \pi(n)=n\}$ to $I_{r,p,n-1}$.

We have:
$$\fix(\pi)=\fix(\pi')+1,$$ $$\exc_A(\pi)=\exc_A(\pi'),$$
$$\csum(\pi)=\csum(\pi'). $$

\item $\pi(n)=n^{\left[\frac{r}{2}\right]}$ and there exists some $k<n$
such that $\pi(k)=k^{\left[\frac{r}{2}\right]}$. Define $\pi'' \in
I_{r,p,n-2}$ as in Section \ref{recurrence}.

Note that $\pi \mapsto \pi''$ is a bijection from the set $\{\pi \in
I_{r,p,n} \mid \pi(n)=n^{\left[\frac{r}{2}\right]}\}$ to
$I_{r,p,n-2}$.

We have:
$$\fix(\pi)=\fix(\pi'')+2,$$ $$\exc_A(\pi)=\exc_A(\pi''),$$
$$\csum(\pi)=\csum(\pi'')+r. $$

\item $\pi(n)=k^{[j]}$ with $k<n$ and we have $\pi(k)=n^{[r-j]}$.
In this case, we use $\pi'' \in I_{r,p,n-2}$ as above. Note that in
this case $\pi \mapsto \pi''$ is a bijection from the set $\{\pi \in
I_{r,p,n} \mid \pi(n)=k^{[j]}\}$ to $I_{r,p,n-2}$. We get in this
case:
$$\fix(\pi)=\fix(\pi''),$$
$$\exc_A(\pi)=\exc_A(\pi'')+\delta_{j,0},$$
$$\csum(\pi)=\csum(\pi'')+r(1-\delta_{j,0}).$$
\end{enumerate}

\medskip

The above consideration gives the following recurrence formula:
\begin{eqnarray*}
f_{r,p,n}(u,v,w) & = & u f_{r,p,n-1}(u,v,w) \\
                    & & +(n-1)(u^2 w^r +(r-1)w^r+v)f_{r,p,n-2}(u,v,w),\quad n\geq1
\end{eqnarray*}

By similar arguments to the ones we have used in Section
\ref{exp_grn}, we get the following generating function and explicit
formula:

\begin{prop}
Let $r\geq1$. The generating function $F_{r,p}(x;u,v,w)$ is given by
$$e^{ux +\frac{1}{2}x^2((u^2+(r-1))w^{r}+v)}$$
\end{prop}

\begin{cor}
The polynomial $f_{r,p,n}(u,v,w)$ is given by
\begin{equation}
\label{pol1}\sum_{j=n/2}^n(n-j)!\binom{n}{n-j,\ n-j ,\
2j-n}\frac{u^{2j-n}(v+(u^2+(r-1))w^r)^{n-j}}{2^{n-j}}.
\end{equation}
\end{cor}

If we substitute $w=1$ and compute the coefficient of $u^mv^{\ell}$
in Formula (\ref{pol1}), we get the following result.

\begin{cor} Let $r\geq1$.  The number of colored involutions in
$G_{r,p,n}$ ($r$ is even, $p \not | \frac{r}{2}$) with exactly $m$
absolute fixed points and ${\rm exc}_A (\pi)=\ell$ is given by
$$ \sumlim_{j=\frac{n}{2}}^{n}\frac{(n-j)!}{2^{n-j}} \binom{n}{n-j,\ 2j-n,\ l,\ \frac{n-m}{2} -l ,\ \frac{m+n}{2}-j}
{(r-1)}^{\frac{n-m}{2}-l} .$$
\end{cor}

For computing the number of involutions with a fixed number of
excedances, we substitute $u=1$ and $v=w^r$ in Formula (\ref{pol1}).

\begin{cor}\label{p dont divide half r}
The number of involutions  $\pi \in G_{r,p,n}$ ($r$ is even, $p \not
| \frac{r}{2}$) with $\exc^{{\rm Clr}}(\pi)=m$ is:
\begin{equation*}
\frac{(\frac{m}{r})!}{2^{\frac{m}{r}}}\binom{n}{\frac{m}{r},\
\frac{m}{r},\ n-2\frac{m}{r}}{(r+1)}^{\frac{m}{r}}
\end{equation*}
\end{cor}


\end{document}